\documentclass[preprint,authoryear,5p,times,twocolumn]{elsarticle}

\usepackage{graphics}
\usepackage{amssymb}
\usepackage{multirow}
\usepackage{color} 
\usepackage{algorithm}
\usepackage{algorithmic}
\usepackage{subcaption}
\usepackage{booktabs}

\begin{document}
\bibliographystyle{model2-names}

\begin{frontmatter}

\title{Heliostat blocking and shadowing efficiency in the video-game era}
\author[label1]{A. Ramos}
\author[label2]{F. Ramos}
\address[label1]{NIC, DESY, Platanenallee 6, 15738 Zeuthen, Germany
    {\tt <alberto.ramos@desy.de>}}
\address[label2]{Nevada Software Informatica S. L. {\tt
    <francisco.ramos@nspoc.com>}}

\begin{abstract}
Blocking and shadowing is one of the key effects in designing and
evaluating a thermal central receiver solar tower plant. Therefore it
is convenient to develop efficient algorithms to compute this
effect. In this paper we explore the possibility of using very
efficient clipping algorithms developed for the video game and imaging
industry to compute the block and shadowing efficiency of a solar
thermal plant layout. We propose an algorithm valid for arbitrary
position, orientation and size of the heliostats. This algorithm
turns out to be very accurate and fast. We show the feasibility of the
use of this algorithm to the optimization of a solar plant by studying
a couple of examples in detail. 
\end{abstract}

\begin{keyword}
optimization; solar thermal electric plant design; field layout; collector
field design.\\
Preprint: {DESY 14-009}
\end{keyword}

\end{frontmatter}

\section{Introduction}

Blocking and shadowing efficiency is one of the main
sources of losses in a field design. Modern plants may have
thousands of heliostats, and this blocking and shadowing efficiency
should be computed for each of the
heliostats. Nevertheless an efficient computation of these effects is non
trivial. It depends not only on the characteristics of the receiver,
but also on the position and orientation of each of the
heliostats. This makes the computation depend on the date and time.

Methods based on ray-tracing techniques achieve a high accuracy on the
estimation of the shadow and block effects, but at a high cost in
computation time. This makes this type of approach unfeasible for
optimization, where we have to evaluate the efficiency of plants with
thousands of heliostats thousands of times. 

A second method to estimate these block and shadow effects is based in
some kind of tessellation of the heliostat (see for
example~\cite{Elsayed:1986,Leonardi:2011}). In this approach the
heliostat is divided in pieces, and then a check is performed to see
if each of these pieces is shadowed or blocked by the rest of the
heliostat field. This approach can in principle achieve an arbitrary
precision by refining the tessellation, but again at the price of a
high computational cost.

A third and practical method used in plant optimization consists in
developing some special formulas for the cases under study, that work
under certain assumptions. Usually one assumes that the heliostat is
blocked and shadowed by other heliostats that are parallel to it, and
neglect the possibility that the same part of the heliostat is
blocked/shadowed by more than one heliostat. While one can argue that
these approximations are accurate enough for the practical purpose of
optimizing a field 
layout, they clearly may have an impact difficult to quantify on the
optimal design. Moreover, modern large plants benefit from 
having more than one tower (see for example~\cite{Ramos:2009aa}), and
therefore some assumptions may be blatantly violated, since
heliostats aiming at different towers are not parallel at all.

It would, therefore, be ideal if one could have an algorithm for evaluating
blocking and shadowing effects that meets the following criteria:
\begin{itemize}
\item Being fast enough that it can be used in optimization process.
\item Being very precise and free of assumptions, working for arbitrary
  positions/orientations of heliostats in a possibly non flat
  ground. 
\end{itemize}

In our opinion none of the previous approaches meet these strict
criteria. Most of the previous approaches to the
problem of evaluating blocking and shadowing efficiencies date back
to the first steps of tower solar plants, around the late 70's
and 80's. But the area of computer graphics has evolved tremendously
since then. Computer video games depend crucially on efficient
\emph{clipping} strategies (determining what portion of a ``picture''
covers some region) to be able to render scenes with high visual
quality. 

From the purely algorithmic point of view, computing the blocking and
shadowing efficiency is nothing more than a clipping problem. More
specifically, we have to determine which portion of the subject heliostat is
visible both from the sun and the receiver point of view. In this work
we apply an efficient 2D polygon clipping algorithm to the
problem of computing the blocking and shadowing heliostat
efficiency. The algorithm will prove to be extremely fast, and works 
for arbitrary shapes, positions, sizes and orientations of the
heliostats (in particular it is as difficult to evaluate for a non
flat ground as for a flat one). It also takes into account correctly
the possibility that part of the blocked/shadowed region may be common
to several heliostats. The only assumption made is that the heliostat
is flat. Since modern heliostats are canted to the slant range
the error made with this approximation is smaller than
1 part in $10^5$ (and, in practice, probably smaller), and therefore
completely negligible. In this way we achieve our goal of a fast,
precise and assumption-free algorithm to evaluate blocking and
shadowing efficiency.

The paper is organized as follows: Section~\ref{sc:gen} contains 
our choice of notation and coordinate systems used for the
computations, as well as some information about polygons, that are a
crucial element in our algorithm. Section~\ref{sc:alg} contains a
detailed description of our algorithm, while section~\ref{sc:test}
contains some examples of our algorithm in use. We finally conclude in
section~\ref{sc:conclusions}.

\section{Generalities}
\label{sc:gen}

\subsection{Coordinate systems}

Vectors are represented with an arrow
$\vec v$, and their components by indices $\vec v = (v_x, v_y, v_z)$. We
will use different coordinate systems. One is the local 
coordinate system of the plant, with the positive $X$ axis
pointing to the north, the positive $Y$ axis pointing to the West and
the $Z$ axis being the height. In this reference system ($S=XYZ$) the
sunlight direction is given by the unit vector  
\begin{equation}
  \label{eq:usun}
  \vec u_s = (-\cos \eta\cos \theta, \cos \eta\sin \theta, -\sin \eta)
\end{equation}
where $\eta$ is the solar height and $\theta$ is the azimuth measured
from the north point. We will consider flat heliostats, and define a
coordinate system for each heliostat, denoted $S'_h=X'Y'Z'$. It has the
origin in the center of the heliostat, and the 
$Z'$ axis being the normal to the heliostat, and the positive
direction in the direction of the reflected light. The $S'_h$
coordinate system is completely determined by saying that the $X'$
axis forms an angle $\phi$ with the $XY$ plane. Obviously the
coordinates of the corners of a flat heliostat 
$h$ in the coordinate system $S_h'$, denoted by $\vec P'^a$ with
$a=1,\dots,N_c$ ($N_c$ is the number of corners of the heliostat), are
given by  
\begin{eqnarray}
  \nonumber
  \vec P'^1 &=& (P'^1_x,P'^1_y,0) \\
  \vec P'^2 &=& (P'^2_x,P'^2_y,0) \\
  \nonumber
  \vdots
\end{eqnarray}
and in the particular case of a rectangular heliostat with dimensions
$L_x\times L_y$ we have
\begin{eqnarray}
  \label{eq:coord}
  \nonumber
  \vec P'^{1} &=& (-L_x/2,L_y/2,0) \\
  \nonumber
  \vec P'^2 &=& (-L_x/2,-L_y/2,0) \\
  \nonumber
  \vec P'^3 &=& (L_x/2,-L_y/2,0) \\
  \vec P'^4 &=& (L_x/2,L_y/2,0)\,.
\end{eqnarray}

If $\vec u_t$ is the unit vector that points \emph{from} the heliostat
\emph{to} the receiver the normal of the heliostat defining the $Z'$
axis is given by
\begin{equation}
  \vec n = \vec u_t - \vec u_s
\end{equation}
that can be used to transform from one coordinate system to the
other. Defining a rotation matrix using the Euler angles in the $ZXZ$
convention  
\begin{eqnarray}
  \nonumber
  R(\alpha,\beta,\gamma) &=& \left(
    \begin{array}{ccc}
      \cos\gamma &\sin\gamma & 0 \\
      -\sin\gamma &\cos\gamma & 0 \\
      0 & 0 & 1 \\
    \end{array}
  \right) \\
  \nonumber
  &\times& \left(
    \begin{array}{ccc}
      1 &0 & 0 \\
      0 & \cos\beta & \sin\beta \\
      0 & -\sin\beta & \cos\beta \\
    \end{array}
  \right) \\
  &\times& \left(
    \begin{array}{ccc}
      \cos\alpha &\sin\alpha & 0 \\
      -\sin\alpha &\cos\alpha & 0 \\
      0 & 0 & 1 \\
    \end{array}
  \right) 
\end{eqnarray}
the application that transform between coordinate systems is given by
\begin{eqnarray}
  \nonumber
  T_h &:& S \longrightarrow S'_h \\
  \vec x' &=& T_h(\vec x) = \tilde R_h\left(\vec n, \phi\right)(\vec
  x-\vec P) 
\end{eqnarray}
where $\vec P$ is the position of the center of the heliostat, $\vec
n$ a vector normal to the heliostat and \footnote{For numerical
  stability the reader is encouraged to 
  implement this rotation matrix either using the $atan2$ function or
  the quaternion form of the rotation matrix.}
\begin{equation}
  \tilde R_h\left(\vec n, \phi\right) = R\left(\arctan\frac{n_1}{-n_2},
    \arctan\frac{\sqrt{n_1^2+n_2^2}}{n_3}, \phi\right).
\end{equation}
It is also straightforward to write the form of the inverse
transformation
\begin{eqnarray}
  \nonumber
  T^{-1}_h &:& S'_h \longrightarrow S \\
  \label{eq:SptoS}
  \vec x &=& T^{-1}_h(\vec x') = \tilde R^T_h\left(\vec n, \phi\right)\vec x' +\vec P 
\end{eqnarray}
that can be used to obtain the coordinates of the corners of an
heliostat in the $S$ coordinate system
\begin{equation}
  \vec P^a = T^{-1}_h(\vec P'^a)\,.
\end{equation}

\subsection{Polygons}

We will make extensive use of simple polygons in the following
discussion. We will represent a polygon by giving the $N_c$ coordinates of
the corners, and use calligraphic letters to identify them. For
example if $\vec P^a$ with $a=1,\dots,4$ are the four corners of an
heliostat, they form the polygon $\mathcal P = \vec P^1\vec P^2\vec
P^3\vec P^4$. We can have polygons in the plane, in which case the
corners are 2D points, or we can have them in space, in which case the
points have three components and it is assumed that all corners lay on
a common plane. 

In the case of 2D polygons, we will write $\vec x \in \mathcal P$ when
the point 
$\vec x$ is inside the polygon $\mathcal P$ and $\vec x \notin
\mathcal P$ otherwise. There is a very efficient
algorithm to determine if a point lies inside a polygon based on
Jordan's curve theorem, called the even-odd rule. One draws a
straight line from the point $\vec x$ to the infinity (in an arbitrary
direction) and count how
many times it crosses the polygon $\mathcal P$. The point is inside
the polygon if and only if this counting is an odd number. We use
simulation of simplicity (SoS) techniques to deal with the possible
degenerate cases~\citep{Mucke:1990}. 

One can define the union, intersection and difference of 2D polygons,
being each of them one or more polygons. The following definitions are
straightforward
\begin{eqnarray}
  \nonumber
  \mathcal A = \mathcal B \cup \mathcal C &\Longleftrightarrow&
  \left\{ 
    \vec x\in \mathcal A \Longrightarrow \vec x \in \mathcal A
    \textrm{ or } \vec x \in \mathcal B
  \right\} \\
  \nonumber
  \mathcal A = \mathcal B \cap \mathcal C &\Longleftrightarrow&
  \left\{ 
    \vec x\in \mathcal A \Longrightarrow \vec x \in \mathcal A
    \textrm{ and } \vec x \in \mathcal B
  \right\} \\
  \mathcal A = \mathcal B \setminus \mathcal C &\Longleftrightarrow&
  \left\{ 
    \vec x\in \mathcal A \Longrightarrow \vec x \in \mathcal A
    \textrm{ and } \vec x \notin \mathcal B
  \right\}.
\end{eqnarray}

\section{Efficient computation of blocking and shadowing
  efficiency and clipping algorithm}  
\label{sc:alg}

\subsection{General idea}

Given a subject heliostat (denoted $c$), we have to determine the
fraction of reflecting surface that is not shadowed or blocked by a
set of heliostats (denoted $h^{(i)}$ with $i=1,\dots,N$). For the sake
of simplicity 
in what follows we will assume that the heliostats are rectangular
with lengths $L_x$ and $L_y$. The generalization to the case of
different lengths for different heliostats or different shapes is
straightforward and will be discussed later.

Our proposed algorithm works in three phases:
\begin{enumerate}
\item Project each of the corners of each heliostat $h^{(i)}$ to the
  subject's heliostat plane. When doing this projection from the sun
  point of view we will determine the shadowing, and when doing it
  from the tower point of view we will determine the blocking. 

  After this first phase we have $2N$ quadrilaterals in the subject's
  heliostat plane.

\item Compute the set-theoretical difference between the subject
  heliostat and the $2N$ quadrilaterals. This is basically a 2D
  polygon clipping problem that can be solved very efficiently. 

  After this second stage we have one polygon (in the general case we
  can have more than one) that represents the effective reflecting
  surface of the heliostat.

\item The shadowing and blocking efficiency is simply the result of
  the ratio between the effective reflective
  surface and the total area of the heliostat.
\end{enumerate}

Now we will describe in detail each of these three steps.

\subsection{Projection to the heliostat plane}
\label{sec:first}

We will consider a square heliostat $c$ of dimensions
$L_x\times L_y$ whose corners have coordinates in the plant
reference system $P^a_c$ with $a=1,2,3,4$. As we have
mentioned the coordinates of these corners in the
reference system ($S_c'$) associated to the heliostat $c$ are given by
Eqs.~\ref{eq:coord}, and these can be transformed to the plant
reference system using Eq.\ref{eq:SptoS}.

To compute the shadowing losses on an heliostat $c$ we have to determine
the intersection of the line that passes through each of the heliostat
corners with direction $\vec u_s$ (given by Eq.~\ref{eq:usun}) with
the  plane of the subject heliostat, defined by having a normal
direction $\vec n_c$ and passing by the center of the heliostat $\vec
X_c$. In this case the intersection points, labeled
$\vec s^a$, are easily determined by
\begin{equation}
  \label{eq:sh}
  \vec s^a = \vec P^a_c + 
  \left(\frac{\vec n_c\cdot \vec X_c - \vec n_c\cdot\vec P^a_c}{\vec
      n_c\cdot \vec u_s} \right)\vec u_s\,.
\end{equation}

In order to compute the blocking losses one has to proceed in a
similar way, but replacing the role of the sun with the role of the
receiver.  In this case we have to compute the intersection of the
lines that pass through the heliostat corners and the receiver at
which each heliostat is aiming (in the general case of a multi-tower
or multi-receiver plant, different heliostats may aim at different
points). Let us denote the position of the receiver in the plant
reference system with $\vec T$. Then, after defining the unit vectors
\begin{equation}
  \vec u^a_t = \frac{\vec T - \vec P^a_c}{|\vec T - \vec P^a_c|}\,,
\end{equation}
the intersection points needed to compute the blocking are given by
\begin{equation}
  \label{eq:bl}
  \vec b^a = \vec P^a_c + 
  \left(\frac{\vec n_c\cdot \vec X_c - \vec n_c\cdot\vec P^a_c}{\vec
      n_c\cdot \vec u_t^a} \right)\vec u_t^a\,.
\end{equation}

The situation where $\vec n_c\cdot \vec u_t^a=0$ (or $\vec n_c\cdot
\vec u_t=0$) represents heliostats laying
perpendicular one to the other, either from the tower (or sun) point
of view. Since in thiss case there is no blocking (or shadowing), one
can discard these situations immediately without any need of computing
any projection. 
\begin{figure}
  \centering
  \includegraphics[width=0.5\textwidth]{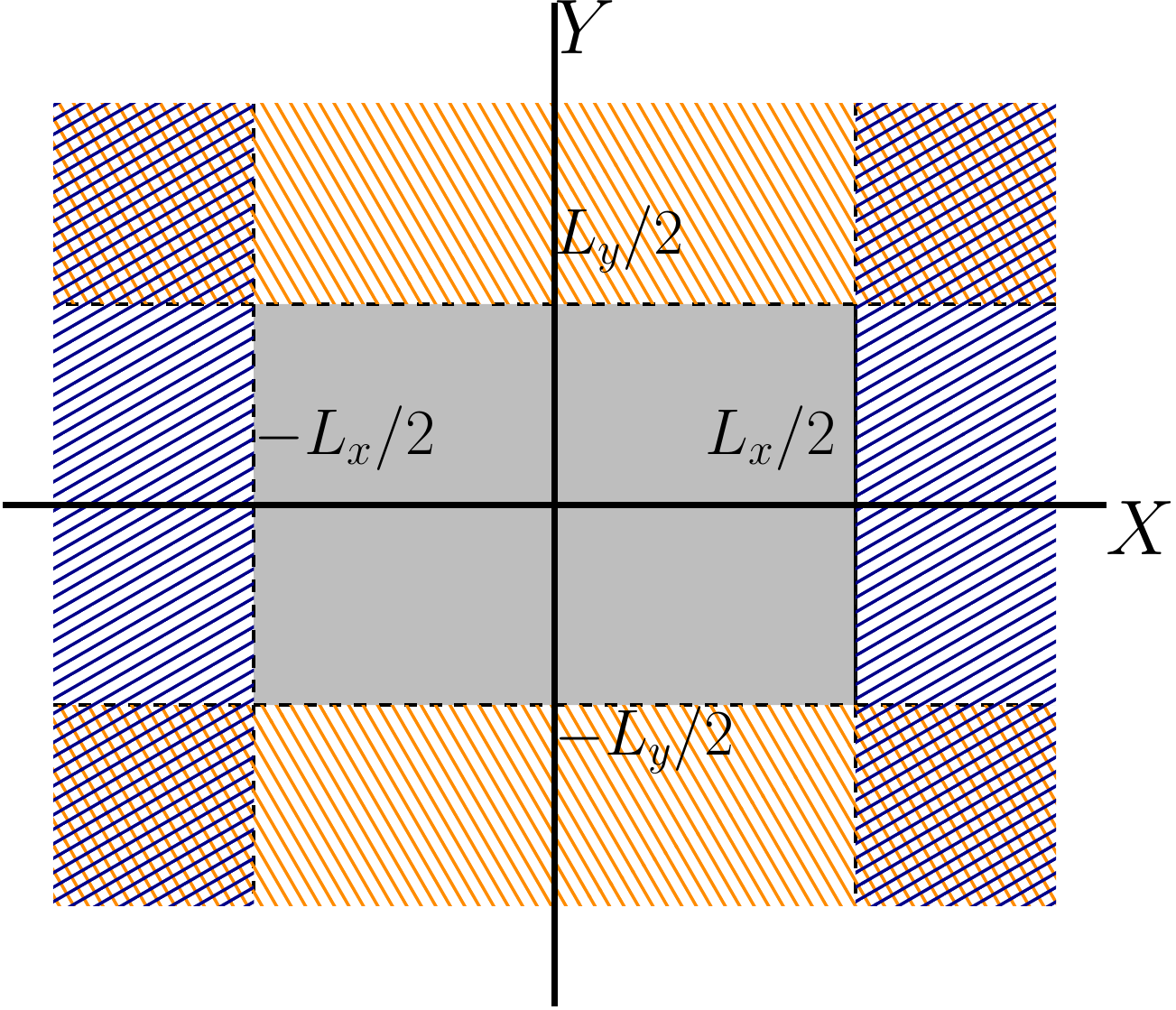}
  \caption{This figure shows the subject heliostat $c$ in the reference
    system defined by itself ($S'_c$). The heliostats that can
    potentially contribute to the losses of the heliostat $c$ define a
  set of quadrilaterals in this reference system via
  Eqs.~\ref{eq:sh} and~\ref{eq:bl}. If any of the these quadrilaterals
  have \emph{all} the corners at the right of the heliostat (i.e. the
  $x$ coordinate of all corners is larger than $L_x/2$), then we know
  that it can not contribute to the heliostat losses. The same
  principle applies when the corners are to the left, top or bottom.} 
  \label{fig:region}
\end{figure}

Each set of points ${\vec s^a}$ and ${\vec b^a}$ determines a
quadrilateral in the heliostat plane, and therefore in the $S'_h$
reference system all these points have a zero value for the third
component. Therefore the first two components of the vectors $\vec
s'^a = T_c(\vec s^a)$ and $\vec b'^a = T_c(\vec b^a)$ define
2D polygons $\mathcal S=\vec s'^1\vec s'^2\vec s'^3\vec s'^4$ and  
$\mathcal B=\vec b'^1\vec b'^2\vec b'^3\vec b'^4$.

This process is repeated for each of the heliostats
$h^{(i)}$, giving a total of $2N$ polygons, that we will denote
$\mathcal S_i$ and $\mathcal B_i$ with $i=1,\dots,N$. 

Note that at this stage we can already eliminate many contributions to
the heliostat losses. If any of the 2D polygons $\mathcal S_i$ or
$\mathcal B_i$ have the four corners in the $S'_c$ reference system in
the same region (labelled $I, II, III$ and $IV$ in
Fig.~\ref{fig:region}), then this polygon can not be the source of
field losses and one can immediately discard it.

\subsection{2D polygon clipping}

Our subject heliostat defines a polygon in his $S'_c$ reference system
given by $\mathcal C = \vec P'^1\vec P'^2\vec P'^3\vec P'^4$. The
total reflecting area $\mathcal R$ is given by 
\begin{equation}
  \mathcal R = \mathcal C \setminus \left(\cup_{i=1}^N
    \mathcal S_i\cup\mathcal B_i
  \right)
\end{equation}
This defines our clipping problem. There are several ways to proceed,
but probably the most effective consists in defining $\mathcal R =
\mathcal C$, and then iteratively subtract to $\mathcal R$ the
polygons $\mathcal B_1,\mathcal S_1,\mathcal B_2,\mathcal
S_2,\dots$ See algorithm~\ref{alg:al}.
\begin{algorithm}
  \centering
  \begin{algorithmic}[1]
    \STATE Compute $\mathcal S_i, \mathcal B_i$ for each heliostat
    $i = 1,\dots, N_h$     \COMMENT{Use Eqs.~\ref{eq:sh} and~\ref{eq:bl}}
    \STATE Start $\mathcal R = \mathcal C$
    \FOR{$n = 1$ to $N_h$}
    \STATE $\mathcal R \leftarrow \mathcal R \setminus \mathcal B_n$
    \STATE $\mathcal R \leftarrow \mathcal R \setminus \mathcal S_n$
    \ENDFOR
    \RETURN $\textrm{Area}(\mathcal R)/\textrm{Area}(\mathcal C)$
    \COMMENT{Use shoelace formula to compute area of polygons.}
  \end{algorithmic}
  \caption{Algorithm to compute block and shadow}
  \label{alg:al} 
\end{algorithm}

To compute efficiently the blocking and shadowing effects it is
crucial to have an efficient algorithm to perform the 2D polygon
clipping (in particular, the polygon set difference). The algorithm
should be valid for arbitrary simple polygons. In principle all the
polygons $\mathcal C$, $\mathcal B_i$ and $\mathcal S_i$ are convex
polygons, but after subtracting some portions of a polygon the result
will in general be a non-convex polygon. In particular we
choose~\citep{Greiner:1998} and the interested reader should consult
the original reference. Being a crucial step in our implementation, we 
will give an overview of how the algorithm works in section~\ref{sc:clip}

\subsection{Final computation}

At this point we have a polygon representing the portion of the
subject heliostat that is visible both from the point of view of the
sun and the point of view of the tower, and this is precisely the
effective reflecting surface. The final efficiency of the
subject heliostat is given by the area of the reflecting
surface over the total area of the heliostat. 

The so called shoelace formula gives the area of an arbitrary simple
2D polygon. If the corners are $\mathcal P = \left\{ \vec P^a;
  a=1,\dots, N_c\right\}$, the area is given by
\begin{equation}
A(\mathcal P) = \frac{1}{2}\sum_{a=1}^N  \left(P_x^aP_y^{a+1} -
  P_x^{a+1}P_y^a \right) \,,
\end{equation}
with $P_x^{N+1} = P_x^1$ and $P_y^{N+1}=P_y^1$. The sign of the area
given by the previous formula gives the orientation of the polygon. If
the points are labeled in counterclockwise direction the sign is
positive and negative otherwise, but the absolute value always
represents the area of the polygon.

A simple application of the previous formula gives the blocking and
shadowing efficiency
\begin{equation}
e = \frac{A(\mathcal R)}{A(\mathcal C)}\,.
\end{equation}

\subsection{Clipping algorithm}
\label{sc:clip}

\begin{figure}
  \centering
  \includegraphics[width=0.3\textwidth]{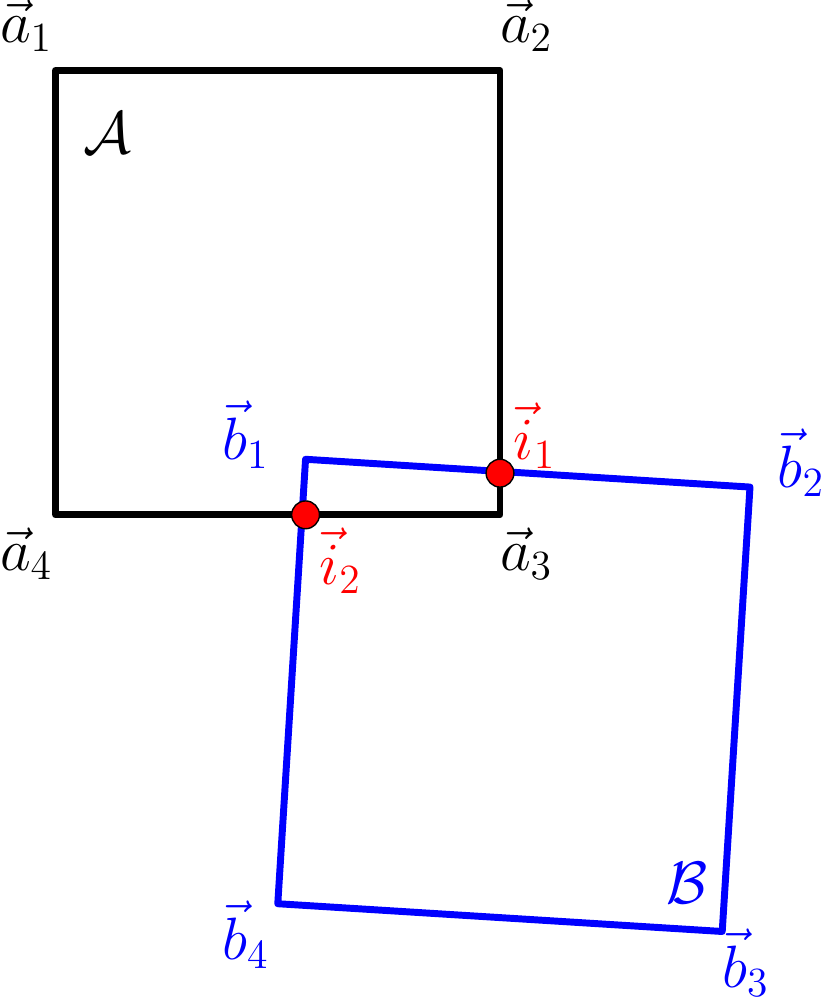}
  \caption{Example of polygon clipping. Given two polygons $\mathcal
    A=\{\vec a_1\vec a_2\vec a_3\vec a_4\}$ and 
    $\mathcal B=\{\vec b_1\vec b_2\vec b_3\vec b_4\}$, to determine
    the difference between the two 
    polygons $\mathcal A\setminus B$ we have first to determine the
    intersection points, labeled $\vec i_1$ and $\vec i_2$. The final
    result is the polygon $\mathcal
A\setminus\mathcal B = \{\vec i_1 \vec b_1\vec i_2\vec a_4\vec a_1\vec a_2\}$.}
  \label{fig:clip}
\end{figure}

The clipping algorithm basically works in three phases. We will refer
to Fig~\ref{fig:clip} as an example, where one
wants to compute the difference between the polygons $\mathcal
    A=\{\vec a_1\vec a_2\vec a_3\vec a_4\}$ and 
    $\mathcal B=\{\vec b_1\vec b_2\vec b_3\vec b_4\}$. 

One starts by determining the intersection points between the two
polygons. In our example there are two of them, and we label them
$\vec i_1$ and $\vec i_2$. These are inserted in
both polygons that now become 
\begin{equation}
\mathcal
A=\{\vec a_1\vec a_2\vec i_1\vec a_3\vec i_2\vec a_4\}
\end{equation}
and 
\begin{equation}
\mathcal B=\{\vec b_1\vec i_1\vec b_2\vec b_3\vec b_4\vec
i_2\}\,.
\end{equation}

The intersection points in each polygon are marked to record
if they are \emph{entry} or \emph{exit} points to the other polygon's
interior. In our example we have
\begin{equation}
\mathcal
A=\{\vec a_1\vec a_2\vec i_1^{\rm entry}\vec a_3\vec i_2^{\rm exit}\vec a_4\}
\end{equation}
and 
\begin{equation}
\mathcal B=\{\vec b_1\vec i_1^{\rm exit}\vec b_2\vec b_3\vec b_4\vec
i_2^{\rm entry}\}\,.
\end{equation}
 
Finally we have to trace the clipped polygon $\mathcal
A\setminus\mathcal B$. The method is quite
simple, one starts in an intersection
point (say for example $\vec i_1$) and moves to the interior of
polygon $\mathcal A$ following polygon $\mathcal B$. When one arrives
to another intersection point, one has to switch polygons: we will
follow polygon $\mathcal A$ moving to the exterior of polygon
$\mathcal B$. In our example this procedure leads to $\mathcal
A\setminus\mathcal B = \{\vec i_1 \vec b_1\vec i_2\vec a_4\vec a_1\vec
a_2\}$, that as the reader can see traces the final polygon
$\mathcal A\setminus\mathcal B$. 

\section{Tests and benchmarks}
\label{sc:test}

\subsection{A simple case}

As a simple test of our method we will consider a solar plant located
at a latitude of $40.08^{\underline{\rm o}}$ with a receiver situated at $100$m of
height. We will compute the blocking and 
shadowing efficiency of an heliostat $c$ at position
\begin{equation}
  c\equiv (108, 0, 5)\,{\rm m}\,.
\end{equation}
This heliostat is surrounded by two heliostats $h_1$ and $h_2$ with
positions
\begin{eqnarray}
  h_1&\equiv& (100, 8, 5)\,{\rm m}\,, \\
  h_2&\equiv& (100, -8, 5)\,{\rm m}\,.
\end{eqnarray}
We will assume that all heliostats have
dimensions $10\,{\rm m}\times 10\,{\rm m}$.

This
example shows nicely how our algorithm can cope with very complex
situations in which the blocking and shadowing effects of several
heliostats overlap. We want to stress that this situation may
happen in real optimal field designs, and that for sure 
this happens very often \emph{during} an optimization process,
where the design might be far from optimal.
\begin{figure}
  \centering
  \includegraphics[width=0.4\textwidth]{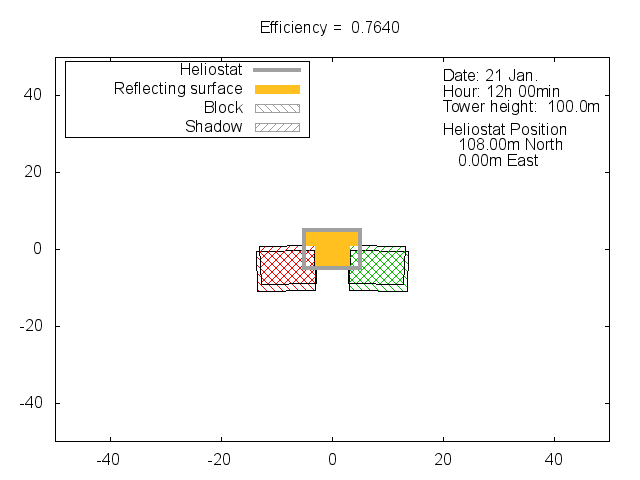}
  \caption{The subject heliostat $c$ and the projection of the shadow
    and block from heliostats $h_1$ (orange) and $h_2$ (blue) at
    noon. The blocking and shadowing efficiency of the heliostat $c$
    is $0.76$. As the reader can see, both heliostats contribute to
    the total blocking and shadowing efficiency of the heliostat $c$. } 
  \label{fig:bs12}
\end{figure}
\begin{figure}
  \centering
  \includegraphics[width=0.4\textwidth]{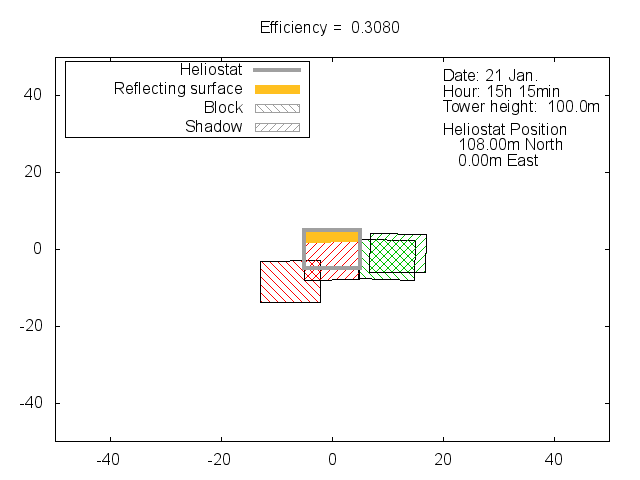}
  \caption{The subject heliostat $c$ and the projection of the shadow
    and block from heliostats $h_1$ (orange) and $h_2$ (blue). The
    blocking and shadowing efficiency of the heliostat $c$ 
    is $0.31$. Both heliostats contribute to
    the total blocking and shadowing efficiency of the heliostat $c$,
    but there is some overlap between the blocking of heliostat $h_1$
    and the shadowing of heliostat $h_2$. } 
  \label{fig:bs15}
\end{figure}

Figure~\ref{fig:bs12} shows the situation at 12h of the 21$^{\rm st}$
of January. As the reader can see both heliostats $h_1$ and $h_2$
contribute to the losses 
by the same amount due to the symmetric position of the heliostats
relative to the receiver and the sun. There is also a large overlap
between the blocked and the shadowed surface of both heliostats. 

Figure~\ref{fig:bs15} shows the same case at 15h 15min of the 21$^{\rm st}$
of January. Again both heliostats $h_1$ and $h_2$ contribute to the
blocking and shadowing efficiency of heliostat $c$. But in this case
the contributions of both heliostats are very different. Most of the
effect comes from the shadowing effects of heliostat $h_2$ while
heliostat $h_1$ only blocks a small portion of heliostat $c$. Moreover
part of this blocking overlaps with the shadowing effect of heliostat
$h_2$. 

As the reader can observe our algorithm deals with these
complicated situations easily and without any particular
constrains. It is as easy to deal with these cases as to
deal with the easy ones in which all the blocking and shadowing
effect comes from a single heliostat. 

\subsection{A real scenario}

\begin{figure}
  \centering
  \includegraphics[width=0.6\textwidth]{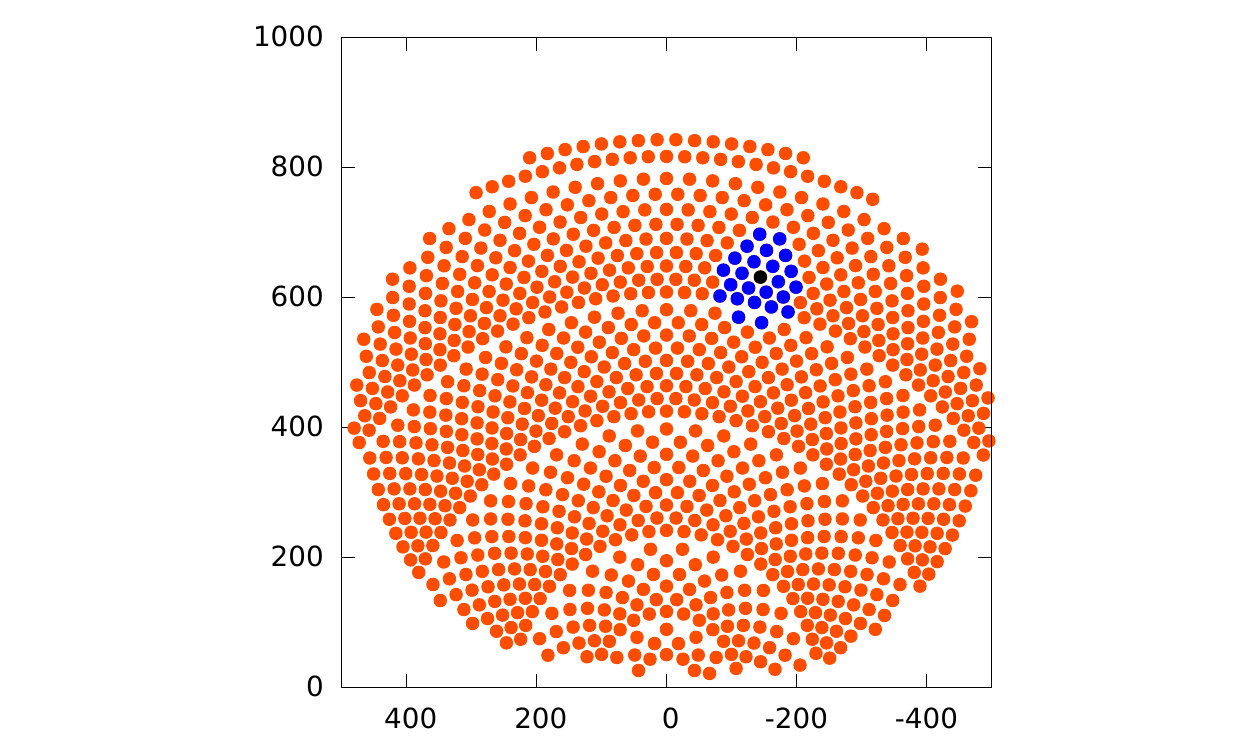}
  \caption{Field layout of a solar plant of 25Mwe with a $150\,m$
    tower. The heliostat marked in black is the subject heliostat that
    we will use to demonstrate our algorithm located at
    $(630.93, -144.41, 5)\,{\rm  m}$. In blue are marked the  
    24 heliostats that are closer to the subject heliostat, and
    therefore the ones that are susceptible of affecting the
    efficiency of the subject heliostat.}
  \label{fig:nspoc}
\end{figure}

As a more realistic example we will consider in detail the computation
of the blocking and shadowing efficiency of one heliostat of a
plant. The field layout was produced along the lines 
of~\citep{Ramos20122536}. An overview of the complete field layout
with the chosen heliostat for the computation of the blocking and
shadowing efficiency can be seen in figure~\ref{fig:nspoc}. In this
example the plant is located at $38.23^{\underline{\rm o}}$ of
latitude and the receiver is located at $150\,{\rm m}$ of height. All
heliostats 
have dimensions $12.88\,{\rm m}\times 9.489\,{\rm m}$.

We will study in detail the blocking and shadowing efficiency of
an heliostat $c$ located at $(630.93, -144.41, 5)\,{\rm  m}$. For
the sake of clarity we will only display the 
24 heliostats (of the 1000 that are in the plant) that are closest
to the subject heliostat. These are the ones that at some moments of
the day affect the blocking and shadowing efficiency and 
are listed in Table~\ref{tab:bs}. We stress that this limitation is
only applied to avoid displaying 999 heliostats in the following figures,
and that the real algorithm, as it is implemented, computes the
effect of all heliostats of the plant. 
\begin{table}
  \centering
  \begin{tabular}{l|l|l|l}
\toprule
n$^{\underline o}$ & Position [m] $(\vec X_c)$ & n$^{\underline o}$ &
Position [m] $(\vec X_c)$\\ 
\midrule
 1  & (614.21, -126.29, 5.00)  & 13  & (660.03, -105.31, 5.00) \\
 2  & (607.84, -153.31, 5.00)  & 14  & (585.22, -161.42, 5.00) \\
 3  & (647.73, -163.37, 5.00)  & 15  & (664.56, -183.31, 5.00) \\
 4  & (654.50, -134.57, 5.00)  & 16  & (678.53, -123.87, 5.00) \\
 5  & (623.77, -172.05, 5.00)  & 17  & (615.45, -199.18, 5.00) \\
 6  & (636.90, -116.27, 5.00)  & 18  & (641.64, -87.72, 5.00) \\
 7  & (591.95, -135.48, 5.00)  & 19  & (689.94, -174.02, 5.00) \\
 8  & (672.18, -153.84, 5.00)  & 20  & (697.15, -143.34, 5.00) \\
 9  & (600.33, -179.84, 5.00)  & 21  & (577.40, -186.87, 5.00) \\
10  & (619.41, -98.828, 5.00)  & 22  & (602.03, -82.30, 5.00) \\
11  & (639.74, -191.65, 5.00)  & 23  & (560.90, -146.37, 5.00) \\
12  & (597.56, -109.09, 5.00)  & 24  & (569.25, -110.83, 5.00) \\
\bottomrule
  \end{tabular}
  \caption{Position of the centers (denoted by $\vec X_c$ in the text)
    of the 24 heliostats of the power plant closer to the
    subject heliostat. These heliostats are marked in blue in
    figure~\ref{fig:nspoc}.}
  \label{tab:bs}
\end{table}

\begin{figure}
  \centering
  \includegraphics[width=0.4\textwidth]{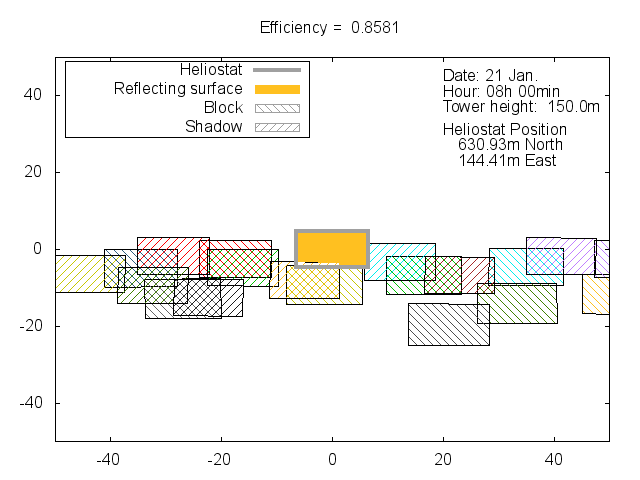}
  \caption{Blocking and shadowing efficiency of an heliostat in a real
    plant situation. At 8h the efficiency is $0.86$. In this partiular
    case two different heliostats produce shadowing losses, and one
    produce some blocking.} 
  \label{fig:bs08}
\end{figure}

\begin{figure}
  \centering
  \includegraphics[width=0.4\textwidth]{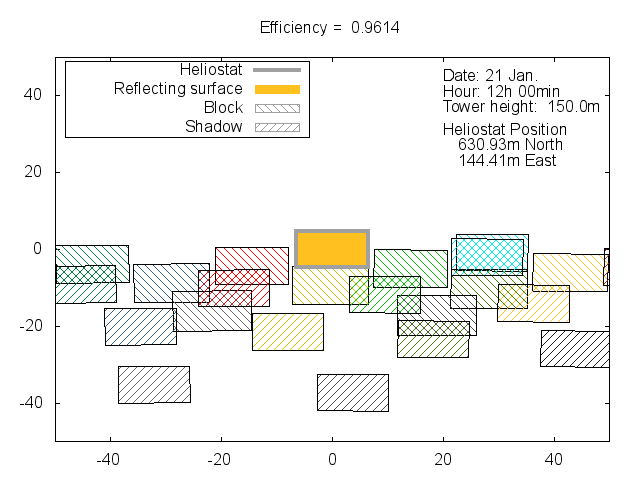}
  \caption{Blocking and shadowing efficiency of an heliostat in a real
    plant situation. At 12h the efficiency is $0.96$, with only one
    heliostat contributing to the losses. } 
  \label{fig:bs122}
\end{figure}

\begin{figure}
  \centering
  \includegraphics[width=0.4\textwidth]{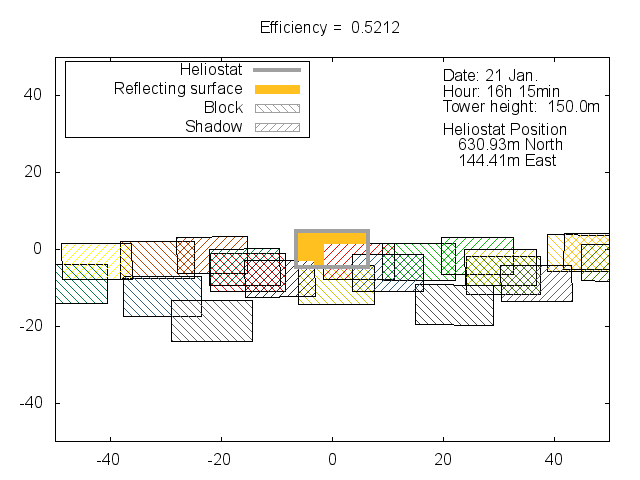}
  \caption{Blocking and shadowing efficiency of an heliostat in a real
    plant situation. At 16h 15min several heliostats contribute to the
    losses and the efficiency drops to $0.52$. This situation is quite
  complex, since there is some overlapping between the blocking
  and shadowing of several heliostats.} 
  \label{fig:bs16}
\end{figure}

We choose the 21$^{\underline{\rm st}}$ of January as the date to
display our results. Figures~\ref{fig:bs08}, \ref{fig:bs122} and
\ref{fig:bs16} 
show the blocking and shadowing at different times of the day. At 8h
(Fig.~\ref{fig:bs08}) in the morning (fig.~\ref{fig:bs08}) the
efficiency is around $0.86$, and two heliostats contribute to the
losses. At 12h (Fig.~\ref{fig:bs122}) only the blocking of one
heliostat contribute to the losses and the efficiency is
$0.96$. Finally at 16h 15m (Fig.~\ref{fig:bs16}) the situation becomes
more complex, since several heliostats contribute to the losses with
some overlap between the blocking and shadowing effects. In this last
case the final efficiency drops to $0.52$. 


From this realistic example we get some important lessons. First, it
is fairly difficult that very many heliostat contribute to the
losses (except just at sunrise or at sunset), and therefore
most of the heliostats of the plant can be easily discarded from the
computation after the first step described in
section~\ref{sec:first}. This step can be implemented regardless of
the algorithm that is used afterwards for the blocking and shadowing
efficiency. Second, it is common, even in an optimized field layout,
that various heliostats contribute to 
the blocking and shadowing efficiency, and that there is some overlap
from the contribution of several heliostats. Third, the heliostats
that contribute to the losses change with the position of the sun. Our
algorithm deals smoothly with these complicated but common 
situations.  

\subsection{Benchmarks}

We have implemented the algorithm previously described in
\texttt{FORTRAN}. For an optimal performance we recommend to use a
linked list to describe polygons in order to avoid moving/copying
data when inserting or removing points from polygons. The original
reference~\citep{Greiner:1998} contains a detailed description of such
data structures (in \texttt{C}). 

As a benchmark we have used the same field layout shown in
figure~\ref{fig:nspoc}. This field consists of 1000 heliostats,
and we have measured the time needed to compute the blocking and
shadowing efficiency of each heliostat. Since this time is very small,
we have repeated the computation 100 times and quote the average as our
final value.
\begin{table}
  \centering
  \begin{tabular}{l|l|l}
    \toprule
    Hour & Execution Time & Field efficiency (Average) \\
    \midrule
    8h 00min & 517 ms & 0.747 \\
    12h 00min & 527 ms & 0.975 \\
    16h 15min & 525 ms & 0.899 \\
    \bottomrule
  \end{tabular}
  \caption{Times for the computation of the blocking and shadowing
    efficiency of each heliostat of the field shown in
    Fig.~\ref{fig:nspoc}.  For each heliostat we have taken into
    account the effect of the remainning 999 heliostats. Since the
    times are very small (around half a second), we have measured the
    time needed to compute 100 times the efficiency of all heliostats
    of the plant, and quote the average as our final result.}
  \label{tab:times}
\end{table}

The results can be seen in Tab.~\ref{tab:times}. The quoted times refer
to runs on a standard laptop (processor \texttt{Intel
  Core(TM) i7-3687U CPU @ 2.10GHz}). In summary we can say 
that, for a plant of a thousand heliostats, our algorithms computes the
blocking and shadowing efficiency in around half a second. These times
seem to be 
pretty independent on the date/time as long as the field layout is
reasonable (i.e. close to optimal, or at least acceptable). 

\section{Conclusions}
\label{sc:conclusions}

In this work we have developed a new algorithm to compute the
heliostat blocking and shadowing efficiency. This algorithm is based on
projecting the image of all heliostats of the field to the plane of
the subject heliostat, both from the point of view of the
sun (shadowing effect) and from that of the tower
(blocking effect). This projected figures form a set of 2D polygons
and our original problem of computing the blocking and shadowing
efficiency is translated into a clipping problem of 2D polygons. This
last problem can be solved very efficiently thanks to the efficient
algorithms developed for the computer video game industry. 

We have used one such algorithm~\citep{Greiner:1998} and showed
that the shadowing and blocking efficiency can be computed very fast:
our \texttt{FORTRAN} code uses less than one second to compute the
blocking and shadowing efficiency of all 1000 heliostats of a plant
running, on a standard laptop. 

The algorithm has no restrictions. Heliostats can be of different
sizes, or be located at different heights. They can aim at different
towers and be located in completely arbitrary positions without any
penalty in the time used by the algorithm to compute the
efficiency. The algorithm is exact, and does not assume any
particular geometry or property of the heliostats. 

We have shown our algorithm at work in some detail for two particular
cases: a simple case in which two heliostats are in front of another,
and a more realistic example in which we have picked up an heliostat
from a real optimal design of a 15MWe plant with 1000 heliostats. We
also want to call the attention  of the reader to two multimedia
files~\cite{mv:bands_564,mv:bands_simple} that show the blocking and
shadowing efficiency for the two cases studied in this paper since the
sunrise till the sunset.

In summary we think that the efficiency of our algorithm makes it a
perfect choice for plant optimization where the blocking and shadowing
efficiency has to be evaluated several thousands of times, while the
lack of any particular assumptions and generality brings peace of mind
to the results of such an optimization process.

\section{Acknowledgments}

The authors want to thank the help of S. Lottini for his comments and
corrections after a careful reading of the manuscript.

\section*{Bibliography}

\bibliography{/home/alberto/docs/bib/math,/home/alberto/docs/bib/campos,/home/alberto/docs/bib/fisica,/home/alberto/docs/bib/computing}

\end{document}